\documentclass{article}
\usepackage{hyperref} 
\usepackage[latin1]{inputenc}  

\usepackage{amsmath,amssymb}

\usepackage{fancyhdr} 
\usepackage{algorithmic,bbm}

\usepackage{graphicx,color} 

\begin{document}
\author
{Alexander Bobenko\thanks{This work is supported by the DFG
  Research Unit \emph{Polyhedral Surfaces.}\hspace{2cm}
   Technische Universität Berlin, Institut für Mathematik,
  Straße des~17.~Juni~136, 10623~Berlin, Germany, \texttt{bobenko,
    schmies@math.tu-berlin.de} } \and Christian
Mercat\thanks{I3M c.c.~51, Université Montpellier 2,
  F-34095 Montpellier cedex 5, France\hspace{2cm}
  \texttt{mercat@math.univ-montp2.fr}} \and Markus
  Schmies$^*$}
\title{Period Matrices of Polyhedral Surfaces} 
\lhead{Period Matrices of Polyhedral Surfaces} 
\rhead{Bobenko, Mercat \& Schmies}
\maketitle
\begin{abstract}
  The linear theory of discrete Riemann surfaces is applied to polyhedral
  surfaces embedded in $\mathbb{R}^3$. As an application we compute
  the period matrices of some classical examples from the surface
  theory, in particular the Wente torus and the Lawson surface.
\end{abstract}
\section{Introduction}
Finding a conformal parameterization for a surface and computing its
period matrix is a classical problem which is useful in a lot of
contexts, from statistical mechanics to computer graphics.

The 2D-Ising model~\cite{M01, CSMcC02, CSMcC03} for example takes
place on a cellular decomposition of a surface whose edges are
decorated by interaction constants, understood as a discrete conformal
structure. In certain configurations, called critical temperature, the
model exhibits conformal invariance properties in the thermodynamical
limit and certain statistical expectations become discrete holomorphic
at the finite level. The computation of the period matrix of higher
genus surfaces built from the rectangular and triangular lattices from
discrete Riemann theory has been addressed in the cited papers by
Costa-Santos and McCoy.

Global conformal parameterization of a surface is important in
computer graphics~\cite{JWYG04, DMA02, BCGB08, TACSD06, KSS06, Wad06} in
issues such as texture mapping of a flat picture onto a curved surface
in $\mathbb{R}^3$. When the texture is recognized by the user as a
natural texture known as featuring round grains, these features should
be preserved when mapped on the surface, mainly because any shear of
circles into ellipses is going to be wrongly interpreted as suggesting
depth increase.  Characterizing a surface by a few numbers is as well
a desired feature in computer graphics, for problems like pattern
recognition. Computing numerically the period matrix of a surface has
been addressed in the cited papers by Gu and Yau.

This paper uses the general framework of discrete Riemann surfaces
theory~\cite{F44,D68,M01,BMS05} and the computation of period matrices
within this framework (based on theorems and not only numerical
analogies). We translate these
theorems into algorithms. Some tests are
performed to check the validity of the approach.

We start with some surfaces with known period matrices and compute
numerically their discrete period matrices, at different level of
refinement. In particular, some genus two surfaces made out of squares and the
Wente torus are considered. We observe numerically good convergence
properties. Moreover, we compute the yet unknown period matrix of the Lawson
surface, recognized it numerically as one of the tested surfaces,
which allowed us to conjecture their conformal equivalence, and
finally prove it.

\section{Discrete conformal structure}

Consider a polyhedral surface in $\mathbb{R}^3$. It has a unique
Delaunay tesselation, generically a triangulation~\cite{BS07}. That is
to say each face is associated with a circumcircle drawn on the
surface and this disk contains no other vertices than the ones on its
boundary. Let's call $\Gamma$ the graph of this cellular
decomposition, $\Gamma_0$ its vertices, $\Gamma_1$ its edges and
complete it into a cellular decomposition with $\Gamma_2$ the set of
triangles.  Each edge $(x,x')=e\in\Gamma_1$ is adjacent to a pair of
triangles, associated with two circumcenters $y,y'$.  The ratio of the
(intrinsic) distances between the circumcenters and the length of the
(orthogonal) edge $e$ is called $\rho(e)$.  \index{circumcenter}

\begin{center}
  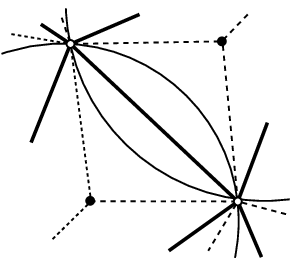
\end{center}

Following~\cite{M01}, we call this data of a graph $\Gamma$, whose
edges are equipped with a positive real number a \emph{discrete
  conformal structure}. A discrete Riemann surface is a conformal
equivalence class of surfaces with the same discrete conformal
structure. It leads to a theory of discrete Riemann surfaces and
discrete analytic functions, developed
in~\cite{F44,D68,M01,M04,M07,BMS05, DKT08}, that shares a lot of
features with the continuous theory and these features are recovered
in a proper refinement limit. We are going to summarize these
results.  \index{discrete conformal structure} \index{conformal
  structure@discrete conformal structure}

In our examples, the extrinsic triangulations are Delaunay. That is to
say the triangulations come from the embedding in $\mathbb{R}^3$ and
the edges $(x,x')\in\Gamma_1$ of the triangulation are the edges of
the polyhedral surface in $\mathbb{R}^3$. On the contrary, the
geodesic connecting the circumcenters $y$ and $y'$ on the surface is
not an interval of a straight line and its length is generically
greater than the distance $||y-y'||$ in $\mathbb{R}^3$. The later
gives a more naive definition of $\rho$, easier to compute that we
will call \emph{extrinsinc}, because contrarily to the
\emph{intrinsic} $\rho$, they depend on the embedding of the surface
in $\mathbb{R}^3$ and are not preserved by isometries. For surfaces
which are refined and flat enough, the difference between extrinsic
and intrinsic distances is not large. We compared numerically the two
ways to compute $\rho$.  The conclusion is that, in the examples we
tested, the intrinsic distance is marginally better, see
Sec.~\ref{sec:WenteExample}.  \index{intrinsic
  distance}\index{distance!intrinsic distance} \index{extrinsic
  distance}\index{distance!extrinsic distance}

The circumcenters and their adjacencies define a 3-valent abstract
(locally planar) graph, dual to the graph of the surface, that we call
$\Gamma^*$, with vertices $\Gamma^*_0=\Gamma_2$, edges
$\Gamma^*_1\simeq \Gamma_1$.  We equip the edge
$(y,y')=e^*\in\Gamma^*_1$, dual to the primal edge $e\in\Gamma_1$,
with the positive real constant $\rho(e^*)=1/\rho(e)$. We define
$\Lambda := \Gamma\oplus\Gamma^*$ the \emph{double} graph, with
vertices $\Lambda_0=\Gamma_0\sqcup\Gamma_0^*$ and edges
$\Lambda_1=\Gamma_1\sqcup\Gamma_1^*$. Each pair of dual edges
$e,e^*\in\Lambda_1$, $e=(x,x')\in\Gamma_1$, $e^*=(y,y')\in\Gamma^*_1$,
are seen as the diagonals of a quadrilateral $(x,y,x',y')$, composing
a quad-graph $\lozenge$, with vertices $\lozenge_0=\Lambda_0$, edges
$\lozenge_1$ composed of couples $(x,y)$ and faces $\lozenge_2$
composed of quadrilaterals $(x,y,x',y')$.  \index{graph!dual graph}
\index{graph!double graph}

The Hodge star, which in the continuous theory is defined by $*(f\, dx
+ g\, dy) = -g\, dx + f\, dy$, is in the discrete case the duality
transformation multiplied by the conformal structure:
\index{Hodge star@discrete Hodge star}
\begin{equation}
  \label{eq:Star}
  \int_{e^*} *\, \alpha :=\rho(e)\int_{e}\alpha
\end{equation}

A $1$-form $\alpha\in C^1(\Lambda)$ is \emph{of type $(1,0)$} if and
only if, for each quadrilateral $(x,y,x',y')\in\lozenge_2$,
$\int_{(y,y')}{\alpha} = i\,\rho(x,x') \int_{(x,x')}{\alpha}$, that is
to say if $*\alpha = -i\, \alpha$. Similarly forms of type
$(0,1)$ are defined by  $*\alpha = +i\, \alpha$. A form is \emph{holomorphic},
resp. anti-holomorphic, if it is closed and of type $(1,0)$, resp. of
type $(0,1)$. A function $f: \Lambda_0\to\mathbb{C}$ is holomorphic
iff $d_\Lambda f$ is.
\index{holomorphic@discrete holomorphic}

We define a wedge product for $1$-forms living whether on edges
$\lozenge_1$ or on their diagonals $\Lambda_1$, as a $2$-form living
on faces $\lozenge_2$. The formula for the latter is:
\begin{align}
  \label{eq:Wedge}
  \iint\limits_{(x,y,x',y')} \alpha\wedge \beta  &:= 
\frac12\left( \int\limits_{(x,x')}\alpha\int\limits_{(y,y')}\beta
-\int\limits_{(y,y')}\alpha\int\limits_{(x,x')}\beta\right)
\end{align}
The exterior derivative $d$ is a derivation for the wedge product, for
functions $f, g$ and a $1$-form $\alpha$:
\begin{align*}
  d(f g) &= f\, dg + g\, df,
& d(f\alpha) = df\wedge\alpha+f\, d\alpha.
\end{align*}
Together with the Hodge star, they give rise, in the compact case, to
the usual scalar product on $1$-forms:
\begin{equation}
  \label{eq:ScalProd}
  \left(\alpha,\,\beta\right):=\iint_{\lozenge_2} \alpha\wedge
  *\bar\beta = 
(*\alpha,\, *\beta) 
= \overline {\left(\beta,\,\alpha\right)}
=\tfrac12\sum_{e\in\Lambda_1}\rho(e)\int_e\alpha\int_e\bar\beta
\end{equation}

The adjoint $d^*=-*\,d\,*$ of the coboundary $d$ allows to define the
discrete Laplacian $\Delta = d^*\,d+d\,d^*$, whose kernel are the
harmonic forms and functions. It reads, for a function at a vertex
$x\in\Lambda_0$ with neighbours $x'\sim x$:
$$\left(\Delta f\right) (x) = \sum_{x'\sim x}\rho(x,x')
\left( f(x) - f(x') \right).$$

\emph{Hodge theorem:} The two $\pm i$-eigenspaces decompose the space
of $1$-forms, especially the space of harmonic forms, into an
orthogonal direct sum.  Types are interchanged by conjugation:
$\alpha\in C^{(1,0)}(\Lambda)\iff \overline\alpha\in
C^{(0,1)}(\Lambda)$ therefore
$$(\alpha,\beta) =  
(\pi_{(1,0)}\alpha,\pi_{(1,0)}\beta)+(\pi_{(0,1)}\alpha,\pi_{(0,1)}\beta)$$
where the projections on $(1,0)$ and $(0,1)$ spaces are
$$\pi_{(1,0)}=\frac 12(\text{Id}+i *),\qquad \pi_{(0,1)}=\frac 12(\text{Id}-i *).$$

The harmonic forms of type $(1,0)$ are the \emph{holomorphic} forms,
the harmonic forms of type $(0,1)$ are the \emph{anti-holomorphic}
forms.

 The $L^2$ norm
of the $1$-form $df$, called the Dirichlet energy of the function $f$,
is the average of the usual Dirichlet energies on each independent
graph \index{Dirichlet energy} \index{energy!discrete Dirichlet energy}
\begin{align}
  E_D(f):=\lVert df\rVert^2 &=
  \left(df,\,df\right)=\frac12\sum_{(x,x')\in\Lambda_1}\rho(x,x')
  \left\lvert f(x') - f(x) \right\rvert^2\label{eq:norm}
  \\
  &=\frac{ E_D(f|_\Gamma)+ E_D(f|_{\Gamma^*})}2 .\notag\end{align} The
conformal energy of a map measures its conformality defect, relating
these two harmonic functions.  A conformal map fulfills the
Cauchy-Riemann equation
\begin{equation}
  \label{eq:CR*}
  *\,df = -i\, df.
\end{equation}
Therefore a quadratic energy whose null functions are the holomorphic
ones is
\begin{equation}
  E_C(f) := \tfrac12\lVert df -i *  df\rVert^2.
\label{eq:EC}
\end{equation}
\index{conformal energy} \index{energy!discrete conformal energy}
It is related to the Dirichlet energy through the same formula as in
the continuous case:
\begin{align}
E_C(f) &= \tfrac12\left( df  -i * df,\, df  -i * df\right)
\notag\\
&= \tfrac12\lVert df \rVert^2+\tfrac12\lVert-i * df \rVert^2
+\, \text{Re}(df,\, -i * df)\notag\\
&=
  \lVert df \rVert^2 + \, \text{Im} \iint_{\lozenge_2} df\wedge\overline{df}
\notag\\
&=
E_D(f) - 2 \mathcal{A}(f)
\label{eq:ECEDA}
\end{align}
where the area of the image of the application $f$ in the complex
plane has the same formulae (the second one meaningful on a simply
connected domain)
\begin{equation}
\mathcal{A}(f) = \frac i2 \iint_{\lozenge_2} df\wedge\overline{df}
= \frac i4 \oint_{\partial\lozenge_2}f\,\overline{df} -\overline{f}\,df
\label{eq:A}
\end{equation}
as in the continuous case.  For a face $(x,y,x',y')\in\lozenge_2$, the
algebraic area of the oriented quadrilateral
$\Bigl(f(x),f(x'),f(y),f(y')\Bigr)$ is given by
\begin{align*}
  \smash{\iint\limits_{(x,y,x',y')}} df\wedge\overline{df}&=
i\, \text{Im}\left(
(f(x')-f(x))\overline{(f(y')-f(y))}
\right)\\
&=-2 i\mathcal{A}\Bigl(f(x),f(x'),f(y),f(y')\Bigr).
\end{align*}

When a holomorphic reference map $z:\Lambda_0\to\mathbb{C}$ is chosen,
a holomorphic (resp. anti-holomorphic) $1$-form $df$ is, locally on
each pair of dual diagonals, proportional to $dz$, resp. $d\bar z$, so
that the decomposition of the exterior derivative into holomorphic and
anti-holomorphic parts yields $df\wedge\overline{df} = \left(|\partial
  f|^2 +|\bar\partial f|^2\right) dz\wedge d\bar z$ where the
derivatives naturally live on faces and are not be confused with the
boundary operator $\partial$.

\section{Algorithm}

The theory described above is straightforward to implement. The most
sensitive part is based on a minimizer procedure which finds the
minimum of the Dirichlet energy for a discrete Riemann surface, given
some boundary conditions. Here is the crude algorithm that we are
going to detail. A normalized homotopy basis  $\aleph$ of a quad-graph
$\lozenge(S)$ is a set of loops $\aleph_k\in \ker \partial$

  \begin{algorithmic}
  \STATE \textbf{Basis of holomorphic forms}{(a discrete Riemann surface $S$)}
  \STATE find a normalized homotopy basis $\aleph$ of $\lozenge(S)$
  \FORALL{$\aleph_k$}
  \STATE compute $\aleph_k^\Gamma$ and $\aleph_k^{\Gamma^*}$    
  \STATE compute the real discrete harmonic form $\omega_k$ on $\Gamma$
    s.t.  $\oint_\gamma \omega_k = \gamma \circ \aleph_k^\Gamma$
    \STATE (check $\omega_k$ is harmonic  on $\Gamma$)
    \STATE compute the form $*\omega_k$ on $\Gamma^*$
    \STATE (check $*\omega_k$ is harmonic  on $\Gamma^*$)
    \STATE compute its holonomies $(\oint_{\aleph_\ell^{\Gamma^*}}
      *\omega_k)_{k,\ell}$ on $\Gamma^*$
    \STATE do likewise for $\omega_k^*$ on $\Gamma^*$
  \ENDFOR  
  \STATE do some linear algebra (R is a rectangular complex matrix) to
    get the basis of holomorphic forms $(\zeta_k)_k = R (\text{Id} + i\,*)
    (\omega_k)_k$ s.t.  $ (\oint_{\aleph_\ell^\Gamma}
    \zeta_k)=\delta_{k,\ell}$
    \STATE define the period matrix
 $\Pi_{k,\ell} := (\oint_{\aleph_\ell^{\Gamma^*}}
      \zeta_k)$
    \STATE do likewise for $(\zeta^*_k)_k$ and  
$\Pi^*_{k,\ell} := (\oint_{\aleph_\ell^{\Gamma}}
      \zeta^*_k)$
  \end{algorithmic}

  Finding a normalized homotopy basis of a connected cellular
  decomposition is performed by several well known algorithms. The way
  we did it is to select a root vertex and grow from there a spanning
  tree, by computing the vertices at combinatorial distance $d$ from
  the root and linking each one of them to a unique vertex at distance
  $d-1$, already in the tree. Repeat until no vertices are left.

  Then we inflate this tree into a polygonal fundamental domain by
  adding faces one by one to the domain, keeping it simply connected:
  We recursively add all the faces which have only one edge not in
  the domain. We stop when all the remaining faces have at least two
  edges not in the domain.

  \begin{center}
    \scalebox{0.6}{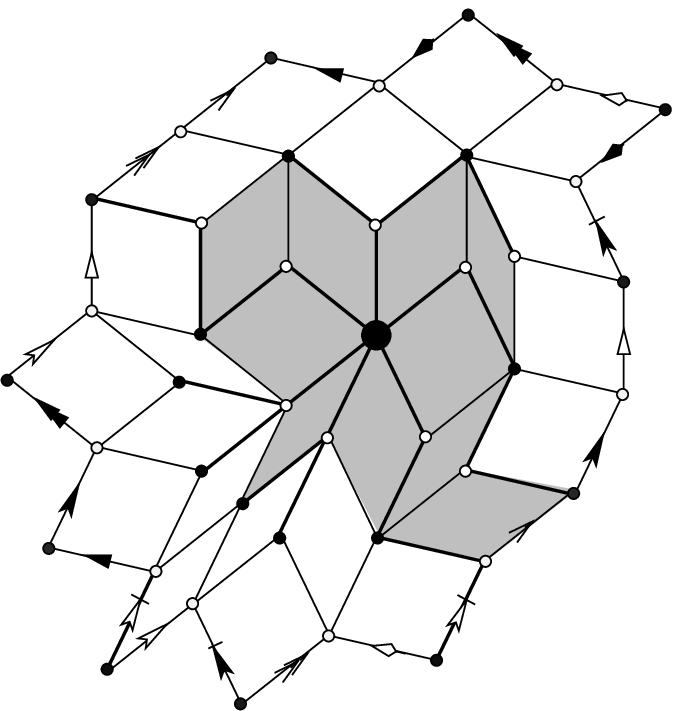}
  \end{center}

  Then we pick one edge (one of the closest to the root) as defining
  the first element of our homotopy basis: adding this edge to the
  fundamental domain yields a non simply connected cellular
  decomposition and the spanning tree gives us a rooted cycle of this
  homotopy type going down the tree to the root. It is (one of) the
  combinatorially shortest in its (rooted) homology class.  We add
  faces recursively in a similar way until we can no go further, we
  then choose a new homotopy basis element, and so on until every face
  is closed. At the end we have a homotopy basis. We compute later on
  the intersection numbers in order to normalize it.
  \begin{center}
    \scalebox{0.4}{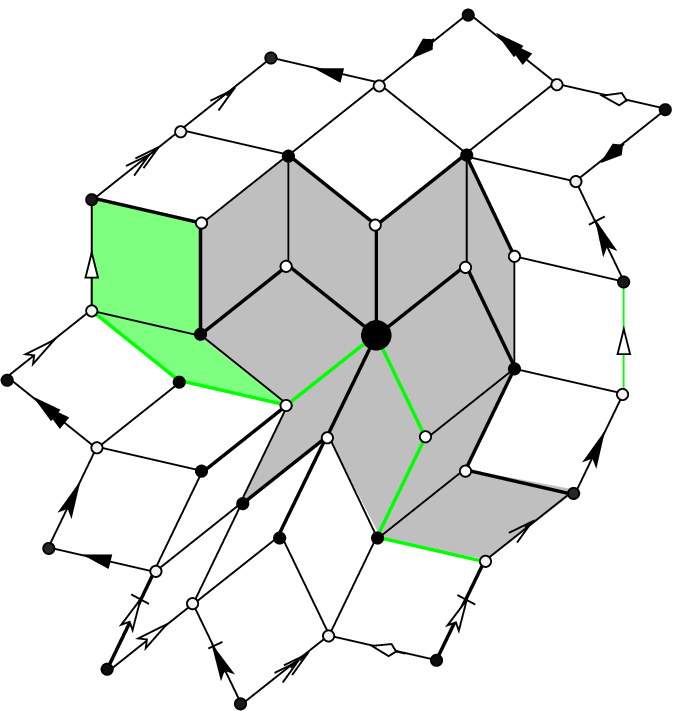
      \hskip 10em  
      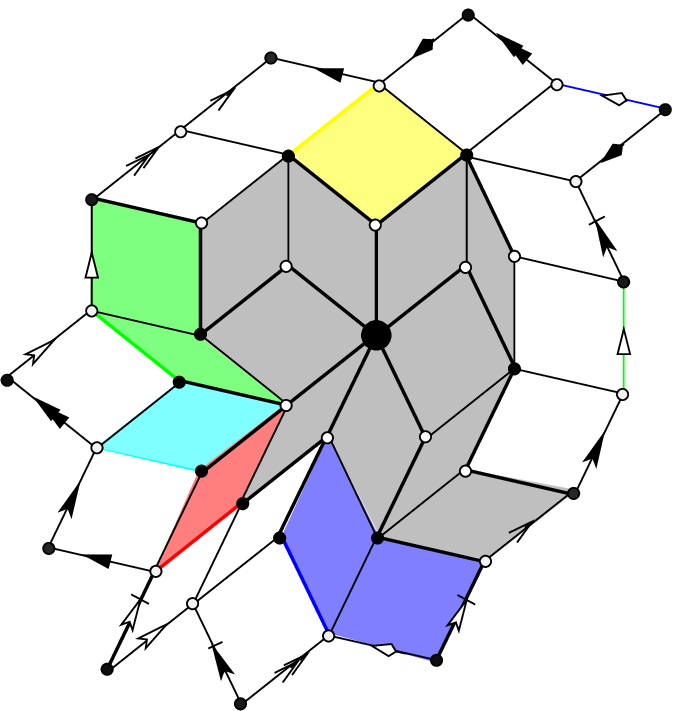
      \hskip 10em
      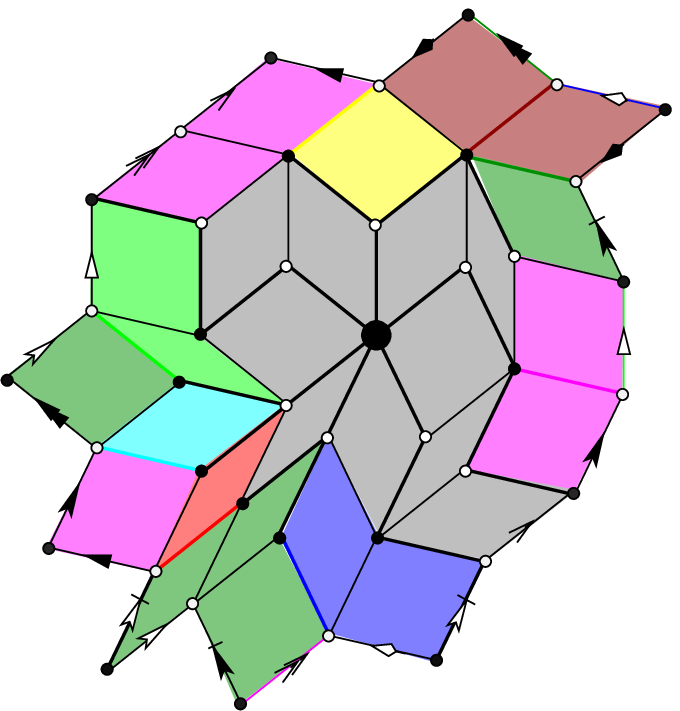}
  \end{center}

  We compute the unique real harmonic form $\eta$ associated with each
  cycle $\aleph$ such that $\oint_\gamma \eta = \gamma \circ \aleph$.
  This is done by a minimizing procedure which finds the unique
  harmonic function $f$ on the graph $\Gamma$, split along $\aleph$,
  whose vertices are duplicated, which is zero at the root and
  increases by one when going across $\aleph$. This is done by linking
  the values at the duplicated vertices, in effect yielding a harmonic
  function on the universal cover of $\Gamma$. The harmonic $1$-form
  $d\, f$ doesn't depend on the chosen root nor on the representative
  $\aleph$ in its homology class.

\section{Numerics}
We begun with testing discrete surfaces of known moduli in order to
investigate the quality of the numerics and the robustness of the
method. We purposely chose to stick with raw \textit{double}
$15$-digits numbers and a linear algebra library which is fast but not
particularly accurate. In order to be able to compare period matrices,
we used a Siegel reduction algorithm~\cite{DHBvHS04} to map them by a
modular transformation to the same canonical form.
\index{Siegel reduction}

\subsection{Surfaces tiled by squares\label{sec:SilholExamples}} Robert
Silhol supplied us with sets of surfaces tiled by squares for which
the period matrices are known~\cite{S06,BS05,S,BS01,RGA97}.  There
are translation and half-translation surfaces: In these surfaces, each
horizontal side is glued to a horizontal side, a vertical to a
vertical, and the identification between edges of the fundamental
polygon are translations for translation surfaces and translations
followed by a half-twist for half-translations. The discrete conformal
structure for these surfaces is very simple: the combinatorics is
given by the gluing conditions and the conformal parameter $\rho\equiv
1$ is constant.

The genus one examples are not interesting because this $1$-form is
then the unique holomorphic form and there is nothing to compute (the
algorithm does give back this known result).  Genus 2 examples are non
trivial because a second holomorphic form has to be computed.

The translation surfaces are particularly adapted because the discrete
$1$-form read off the picture is already a discrete holomorphic form.
Therefore the computations are accurate even for a small number of
squares. Finer squares only blur the result with numerical noise. For
half-translation surfaces it is not the case, a continuous limit has
to be taken in order to get a better approximation.

\begin{center}
  \begin{tabular}{|c|c|}
    \hline
    Surface \& Period Matrix  & Numerical Analysis\\
    \hline
    $\vcenter{\hbox to 11em {\hfill\resizebox{!}{5em}{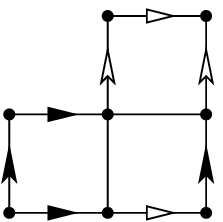}\hfill}}$
     $\Omega_1 = \frac{i}{3} \left(
    \begin{array}{cc}
      5 & -4\\
      -4 & 5
    \end{array} \right)$ & \begin{tabular}{c|c}
     \#vertices & $\left\| \Omega_D - \Omega_1 \right\|_{\infty}$\\
      \hline
      25 & $1.13 \cdot 10^{^{- 8}}$\\
      \hline
      106 & $3.38 \cdot 10^{^{- 8}}$\\
      \hline
      430 & $4.75 \cdot 10^{^{- 8}}$\\
      \hline
      1726 & $1.42 \cdot 10^{^{- 7}}$\\
      \hline
      6928 & $1.35 \cdot 10^{^{- 6}}$
    \end{tabular}\\
    \hline
    $\vcenter{\hbox to 9em {\hfill\resizebox{!}{6em}{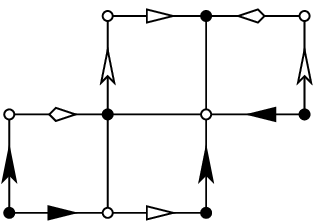}\hfill}}$
  $\Omega_2 = \frac{1}{3} \left(
    \begin{array}{cc}
      - 2 + \sqrt{8} i & 1 - \sqrt{2} i\\
      1 - \sqrt{2} i & -2 + \sqrt{8} i
    \end{array} \right)$ & \begin{tabular}{c|c}
      \hline
      \#vertices & $\left\| \Omega_D - \Omega_2 \right\|_{\infty}$\\
      \hline
      14 & $3.40 \cdot 10^{^{- 2}}$\\
      \hline
      62 & $9.51 \cdot 10^{^{- 3}}$\\
      \hline
      254 & $2.44 \cdot 10^{^{- 3}}$\\
      \hline
      1022 & $6.12 \cdot 10^{^{- 4}}$\\
      \hline
      4096 & $1.53 \cdot 10^{^{- 4}}$
    \end{tabular}\\
    \hline
    $\vcenter{\hbox to 10em {\hfill\resizebox{!}{7.5em}{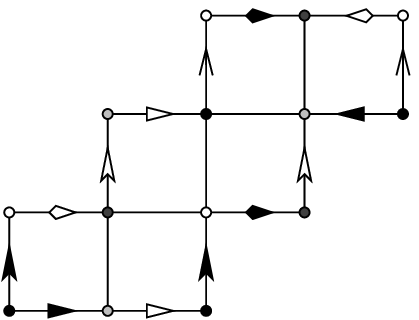}\hfill}}$
   $\Omega_3 = \frac{i}{\sqrt{3}} \left(
    \begin{array}{cc}
      2  &  -1\\
      -1 &  2
    \end{array} \right)$ & \begin{tabular}{c|c}
      \#vertices & $\left\| \Omega_D - \Omega_3 \right\|_{\infty}$\\
      \hline
      22 & $3.40 \cdot 10^{^{- 3}}$\\
      \hline
      94 & $9.51 \cdot 10^{^{- 3}}$\\
      \hline
      382 & $2.44 \cdot 10^{^{- 4}}$\\
      \hline
      1534 & $6.12 \cdot 10^{^{- 5}}$\\
      \hline
      6142 & $1.53 \cdot 10^{^{- 6}}$
    \end{tabular}\\
    \hline
  \end{tabular}
\end{center}

Using $15$ digits numbers, the theoretical numerical accuracy is
limited to $8$ digits because our energy is quadratic therefore half
of the digits are lost. Using an arbitrary precision toolbox or
Cholesky decomposition in order to solve the linear system would allow
for better results but it is not the point here.

\subsection{Wente torus\label{sec:WenteExample}}
\index{Wente torus}
\begin{figure}[htbp]
  $\begin{array}{cc}
    \includegraphics[width=0.5\columnwidth]{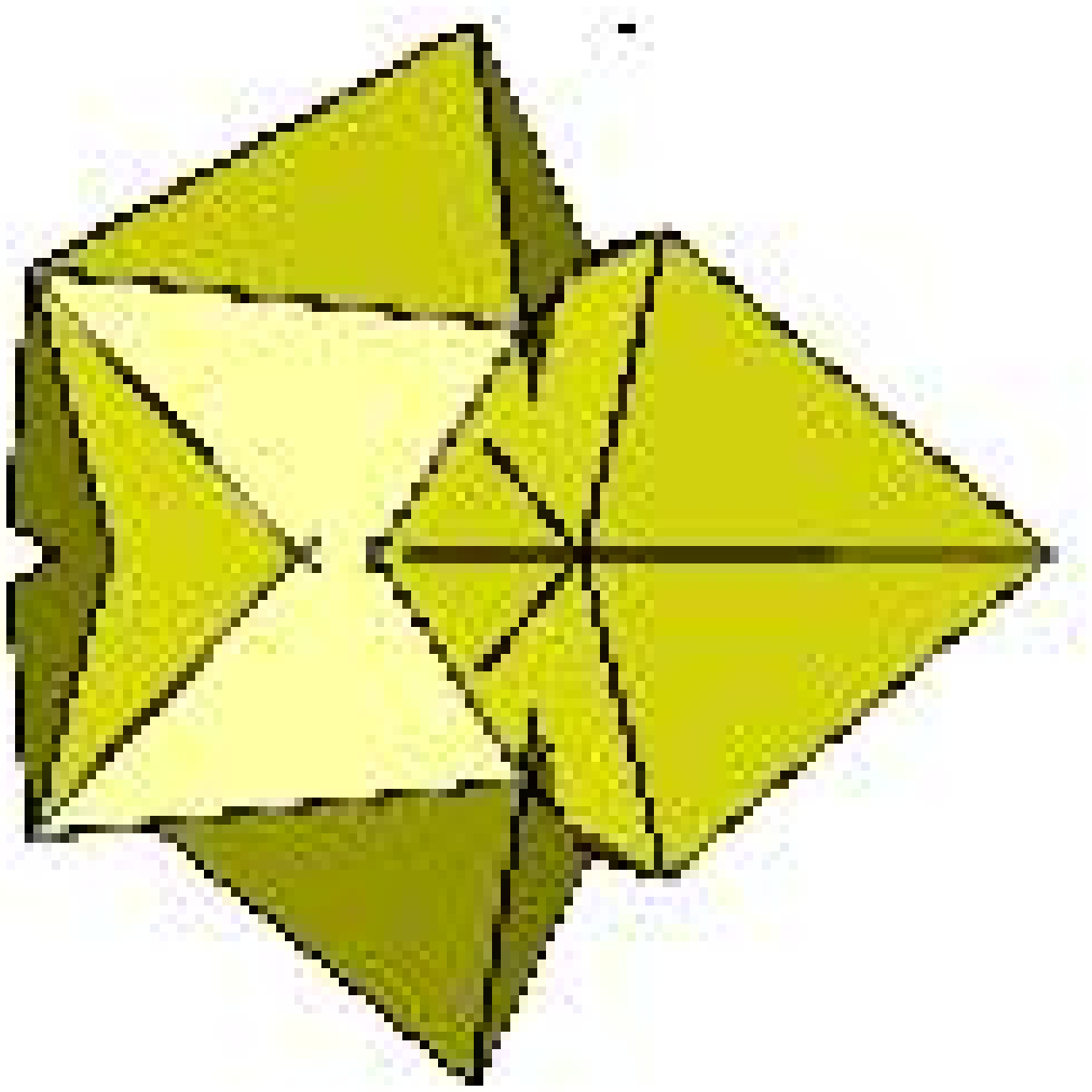} &
    \includegraphics[width=0.5\columnwidth]{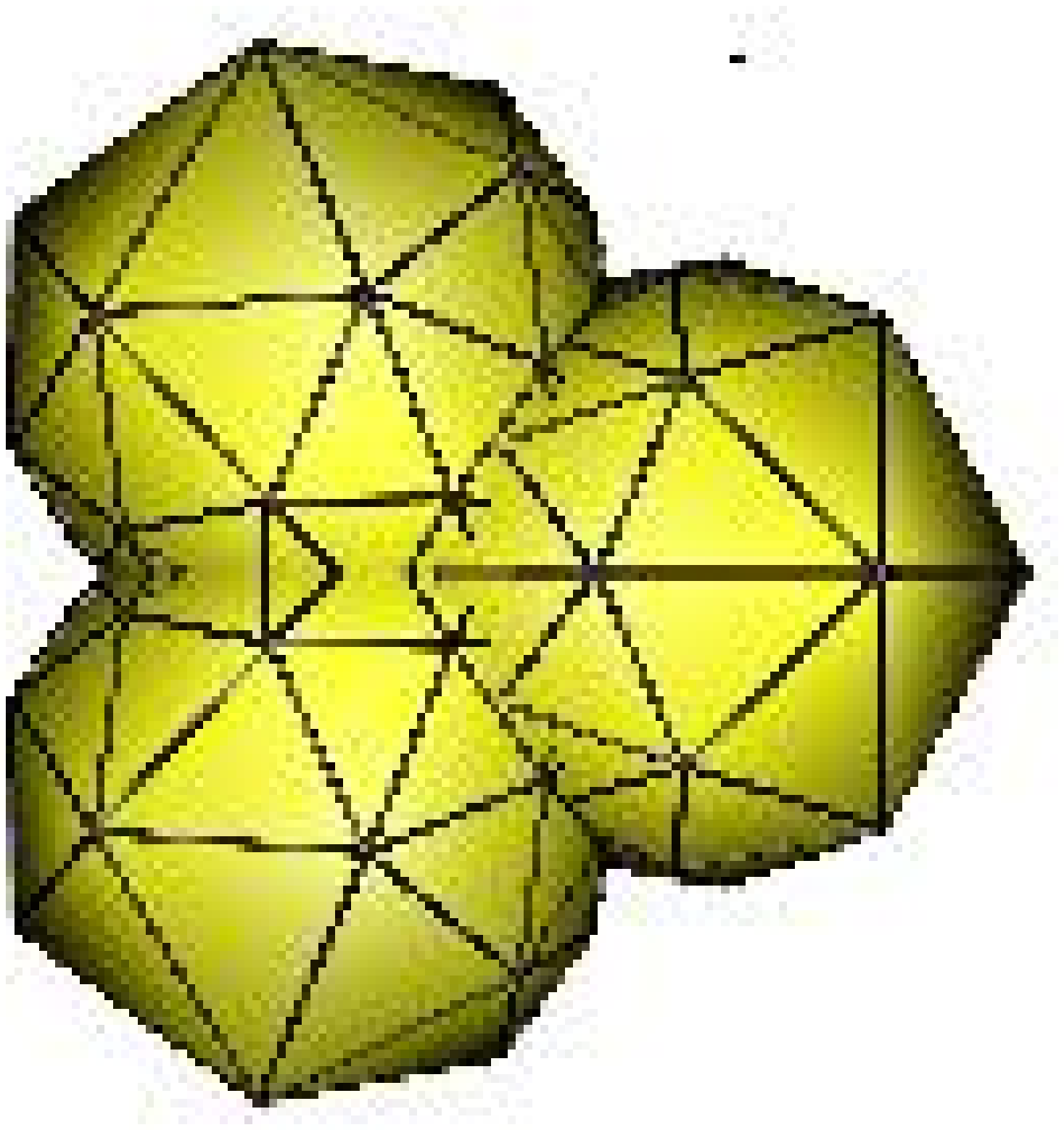}\\
      
        \text{Grid} : 10 \times 10
     
 & 
      \text{Grid} : 20 \times 20
    \\
    \includegraphics[width=0.5\columnwidth]{pictures/wente400Top.eps} &
    \includegraphics[width=0.5\columnwidth]{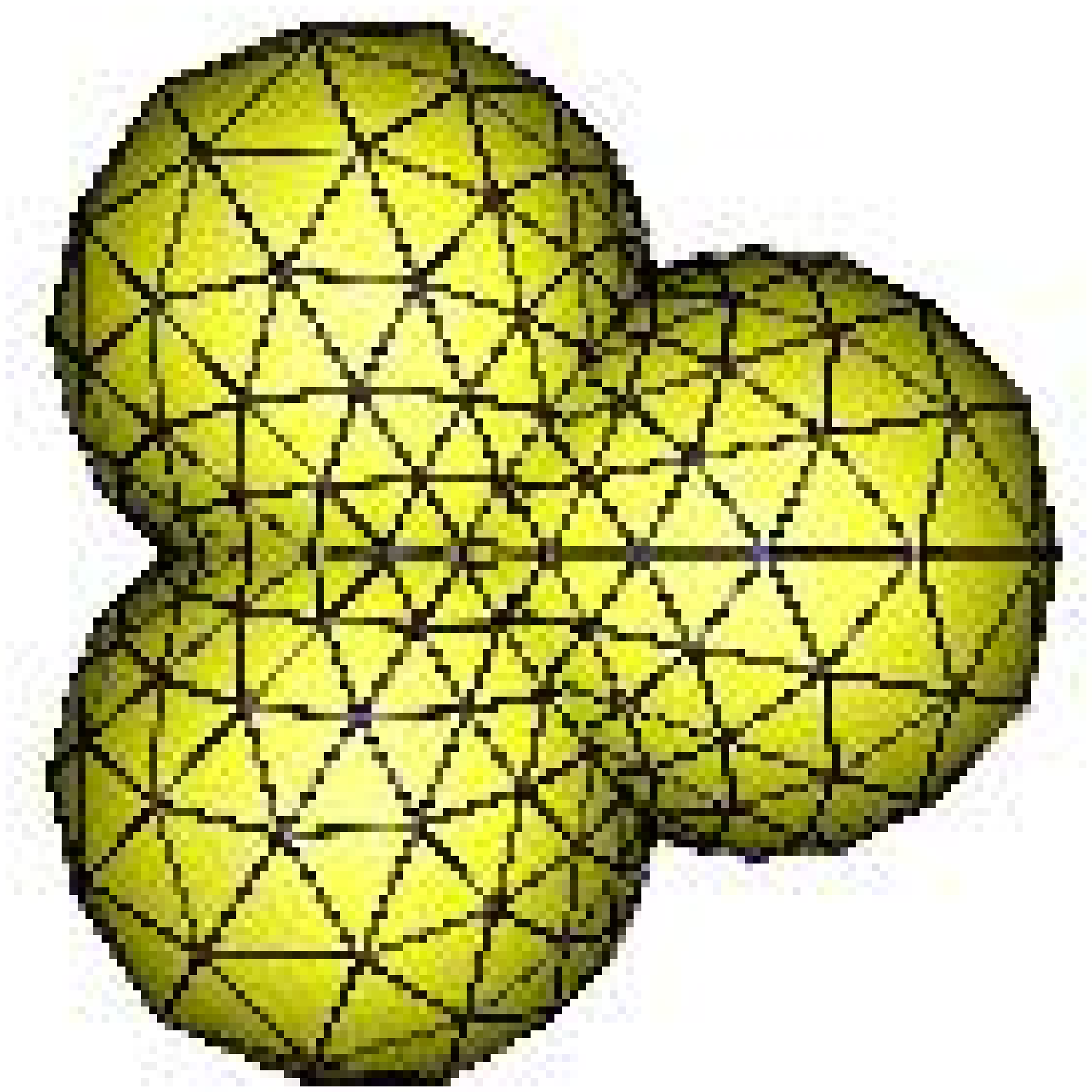}\\
    
      \text{Grid} : 40 \times 40
     &  
      \text{Grid} : 80 \times 80
    
  \end{array}$
  \caption{\label{fig:WenteTori}Regular Delaunay triangulations of the Wente
  torus}
\end{figure}

For a first test of the numerics on a an immersed surface in
$\mathbbm{R}^3$ our choice is the famous CMC-torus discovered by
Wente~\cite{W86} for which an explicit immersion formula exists in terms
of theta functions~\cite{B91}. The modulus of the rhombic Wente torus
can be read from the immersion formula:

\[ \tau_w \approx 0.41300 \ldots + 0.91073 \ldots i \approx \exp (i
1.145045 \ldots .) . \] 

We compute several regular discretization of the Wente torus
(Fig.~\ref{fig:WenteTori}) and generate discrete conformal structures
using $\rho_{\text{ex}}$ that are imposed by the extrinsic Euclidean
metric of $\mathbbm{R}^3$ as well as \ $\rho_{\text{in}}$ which are
given by the intrinsic flat metric of the surface.  For a sequence of
finer discretizations of a smooth immersion, the two sets of numbers
come closer and closer. For these discrete conformal structures we
compute again the moduli which we denote by $\tau_{\text{ex}}$ and
$\tau_{\text{im}}$ and compare them with $\tau_w$ from above:

\begin{center}
  \begin{tabular}{|c|c|c|}
    \hline
    ~~~ Grid ~~~ & ~~~ $\left\| \tau_{\text{in}} - \tau_w \right\|$
    ~~~ & ~~~
    $\left\| \tau_{\text{ex}} - \tau_w \right\|$ ~~~ \\
    \hline
    $10\times 10$ & $5.69 \cdot 10^{- 3} $ & $5.00 \cdot 10^{- 3}$\\
    \hline
    $20\times 20$ & $2.00 \cdot 10^{- 3}$ & $5.93 \cdot 10^{- 3}$\\
    \hline
    $40\times 40$ & $5.11 \cdot 10^{- 4}$ & $1.85 \cdot 10^{- 3}$\\
    \hline
    $80\times 80$ & $2.41 \cdot 10^{- 4}$ & $6.00 \cdot 10^{- 4}$\\
    \hline
  \end{tabular}
\end{center}

For the lowest resolution the accuracy of $\tau_{\text{ex}}$ is slightly
better then the one of $\tau_{\text{in}}$. For all other the discrete
conformal structures with the intrinsically generated $\rho_{\text{in}}$ yields
significant higher accuracy.

\subsection{Lawson surface}

\begin{figure}[htbp]
  $\begin{array}{cc}
    \scalebox{0.75}{\includegraphics[width=0.5\columnwidth]
                                    {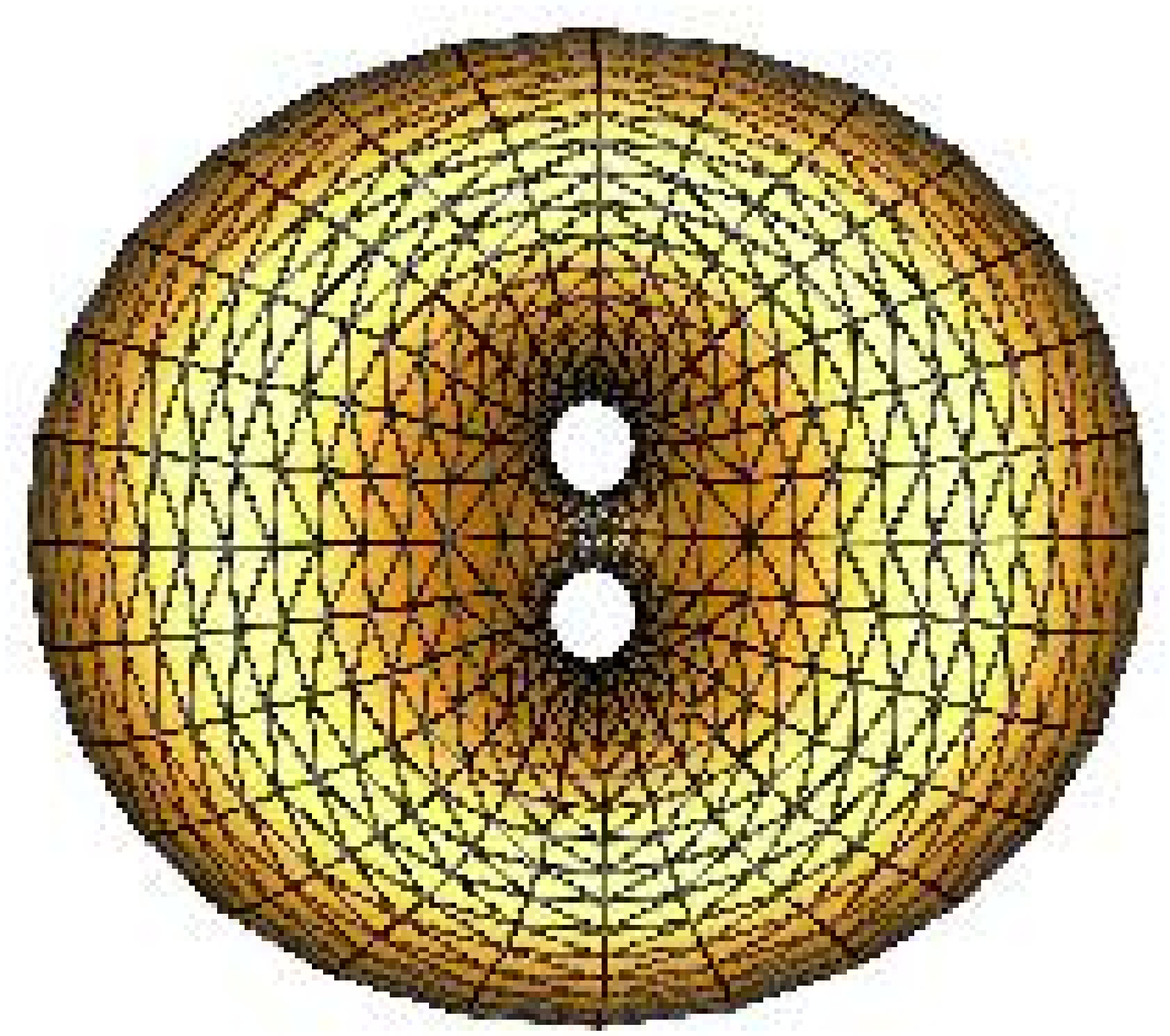}} &
    \scalebox{0.75}{\includegraphics[width=0.5\columnwidth]
                                    {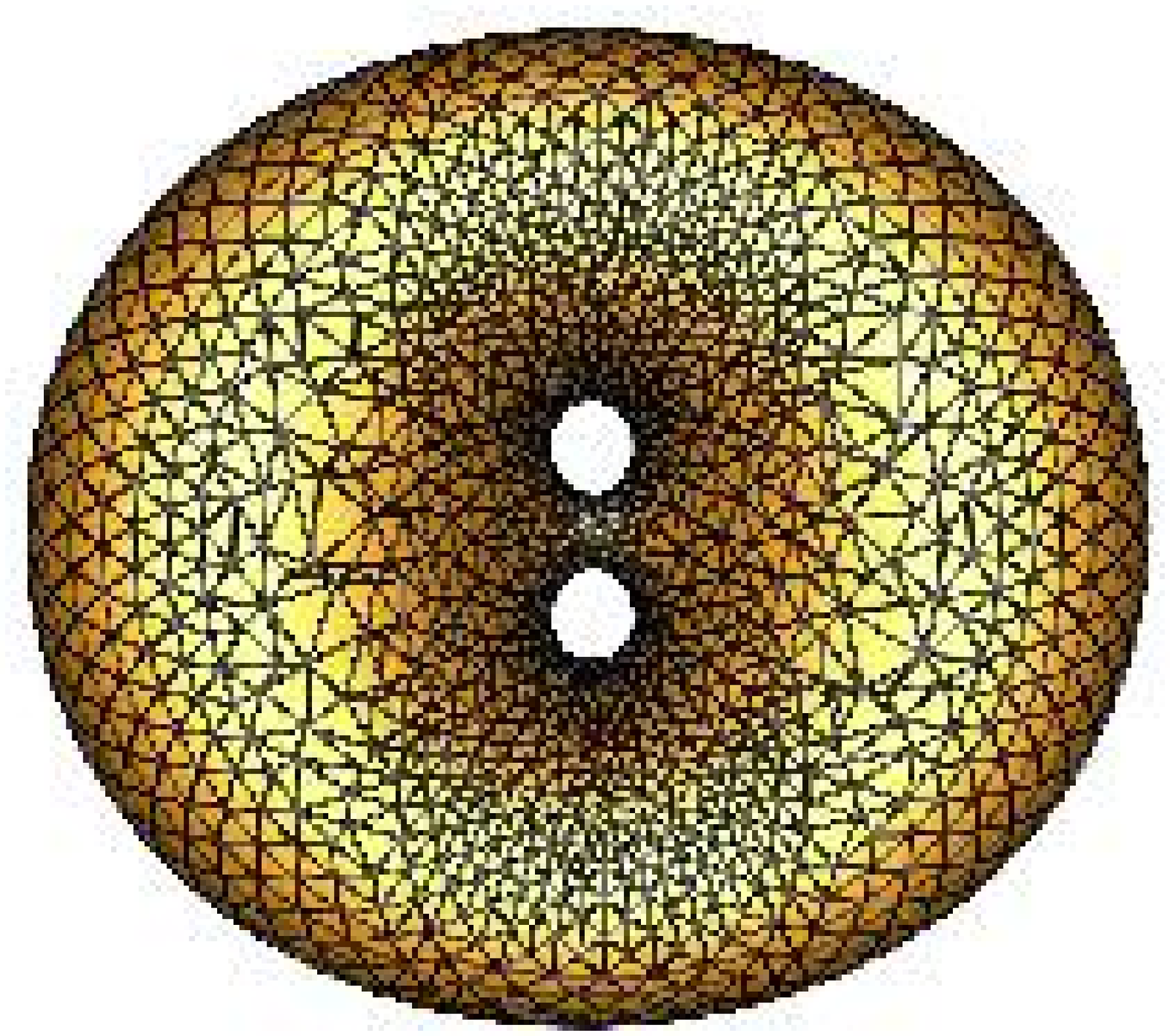}}\\
      \text{1162 vertices}
    & 
      \text{2498 vertices}
    \\
    & \\
    & 
  \end{array}$
  \caption{\label{fig:LawsonSurface}Delaunay triangulations of the Lawson surface}
\end{figure}
\begin{figure}[hbtp]
   \centering
\scalebox{0.7}{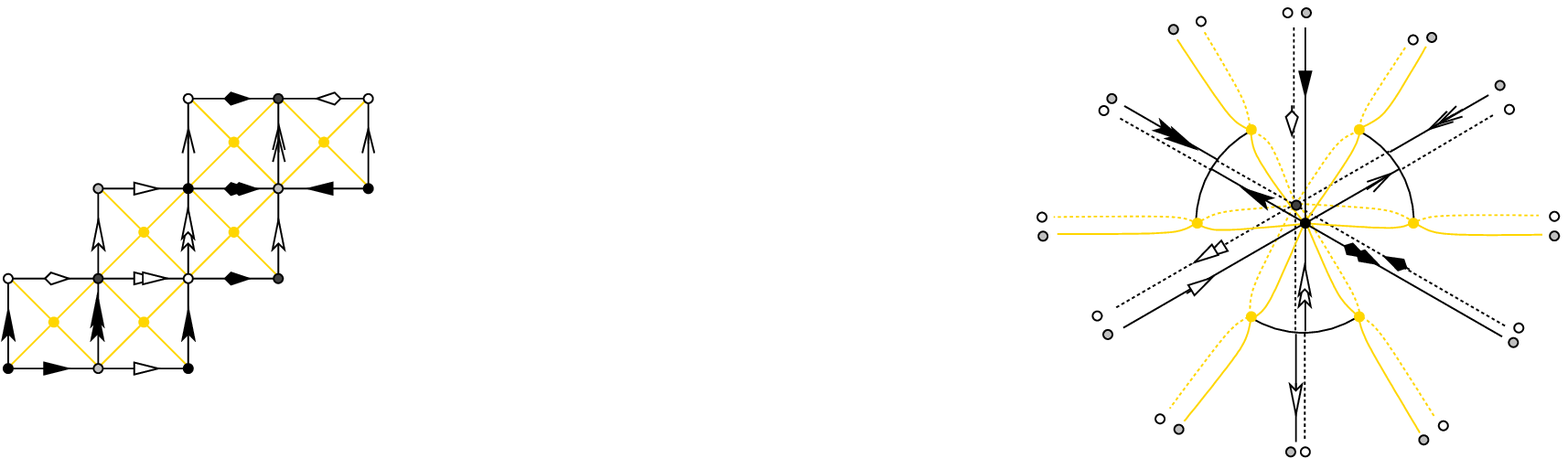}   
\scalebox{2.3}{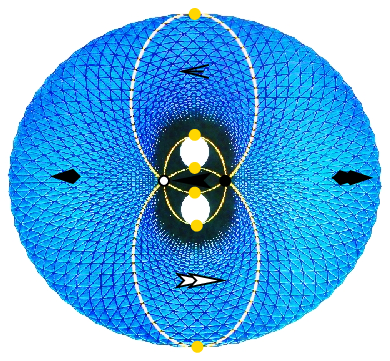}  
\scalebox{2.3}{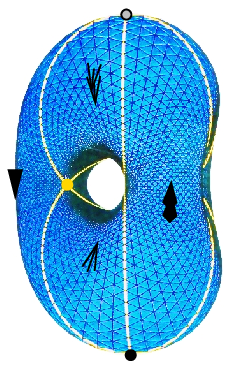}
   \caption{The Lawson surface is conformally equivalent to a surface made of squares.}
   \label{fig:LawsonQuad}
 \end{figure}

\index{Lawson surface} Finally we apply our method to compute the
period matrix of Lawson's  genus 2 minimal
surface in $\mathbb{S}^3$~\cite{Law70}. Konrad Polthier~\cite{PP93} supplied us with
several resolution of the surface which are generated by a coarsening
and mesh beautifying process of a very fine approximation of the
Lawson surface (Fig.~\ref{fig:LawsonSurface}).  Our numerical analysis
gives evidence that the period matrix of the Lawson surface is
\[ \Omega_l = \frac{i}{\sqrt{3}} \left( \begin{array}{cc}
     2 & - 1\\
     - 1 & 2
   \end{array} \right) \]
 which equals the period matrix $\Omega_3 $ of the third example from
 Sec.~\ref{sec:SilholExamples}. Once conjectured that these two
 surfaces are conformally equivalent, it is a matter of checking that
 the symmetry group of the Lawson genus two surface yields indeed this
 period matrix, which was done, without prior connection,
 in~\cite{BBM85}. An explicit conformal mapping of the surfaces can be
 found manually: The genus 2 Lawson surface exhibits by construction
 four points with an order six symmetry and six  points of order four,
 which decomposes the surface into 24 conformally equivalent
 triangles, of angles $\frac{\pi}6$,  $\frac{\pi}2$,  $\frac{\pi}2$.
 Therefore an algebraic equation for the Lawson surface is
 $y^2=x^6-1$, with six branch points at the roots of unity. The
 correspondance between the points in the square picture of the
 surface and the double sheeted cover of the complex plane is done in
 Fig.~\ref{fig:LawsonQuad}. In particular the center of the six
 squares are sent to the branch points, the vertices are sent to the
 two copies of $0$ (black and dark gray) and $\infty$ (white and light
 gray), the square are sent to double sheeted two gons corresponding
 to a sextant.

 Similarly to Sec.~\ref{sec:WenteExample} we compute the period
 matrices $\Omega_{\text{ex}}$ and $\Omega_{\text{in}}$ for different
 resolutions utilizing weights imposed by the extrinsic and intrinsic
 metric and compare the results with our conjectured period matrix for
 the Lawson surface $\Omega_l$:

\begin{center}
  \begin{tabular}{|c|c|c|}
    \hline
    ~~~ \#\text{vertices} ~~~  & ~~~  $\left\| \Omega_{\text{in}} - \Omega_l \right\|_{\infty}$  ~~~ & ~~~ 
    $\left\| \Omega_{\text{ex}} - \Omega_l \right\|_{\infty}$ ~~~ \\
    \hline
    $1162$ & $1.68 \cdot 10^{- 3} $ & $1.68 \cdot 10^{- 3}$\\
    \hline
    $2498$ & $3.01 \cdot 10^{- 3}$ & $3.20 \cdot 10^{- 3}$\\
    \hline
    $100$90 & $8.55 \cdot 10^{- 3}$ & $8.56 \cdot 10^{- 3}$\\
    \hline
  \end{tabular}
\end{center}

Our first observation is that the matrices $\Omega_{\text{ex}}$ and
$\Omega_{\text{in}}$ almost coincide. Hence the method for computing
the $\rho$ seems to have only little influence on this result (compare
also Sec.~\ref{sec:WenteExample}). Further we see that figures of the
higher resolution surface, i.e. with 2498 and 10090 vertices ar worse
than the coarsest one with 1162 vertices. The mesh beautifying process
was most successful on the coarsest triangulation of the Lawson
surface (Fig.~\ref{fig:LawsonSurface}). The quality of the mesh has a
significant impact on the accuracy of our computation: One can see
that the triangles on the coarsest example are of even shapes with
comparable side lengths, while the finer resolution contains thin
triangles with small angles. The convergence speed proven
in~\cite{M01} is governed by this smallest angles, accounting for the
poor result. Therefore for this method to be applicable, the data
should be well suited and it is not enough to have very refined data
if the triangles themselves are not of a good shape. A further
impediment to the method is the fact that the triangulation should be
Delaunay. If it is not, it can be repaired by the algorithm described
in~\cite{FSSB07}.

\bibliographystyle{abbrv}
 \bibliography{surface}

\begin{thebibliography}{10}

\bibitem{BBM85}
M.~V. Babich, A.~I. Bobenko, and V.~B. Matveev.
\newblock Solution of nonlinear equations, integrable by the inverse problem
  method, in {J}acobi theta-functions and the symmetry of algebraic curves.
\newblock {\em Izv. Akad. Nauk SSSR Ser. Mat.}, 49(3):511--529, 1985.

\bibitem{BCGB08}
M.~Ben-Chen, C.~Gotsman, and G.~Bunin.
\newblock
  \href{http://www.cs.technion.ac.il/~gotsman/AmendedPubl/Miri/EG08_Conf.pdf}{%
Conformal Flattening} by curvature prescription and metric scaling.
\newblock In {\em Computer Graphics Forum}, 27(2). Proc. Eurographics, 2008.

\bibitem{B91}
A.~I. Bobenko.
\newblock All constant mean curvature tori in {${\bf R}\sp 3,\;S\sp 3,\;H\sp
  3$} in terms of theta-functions.
\newblock {\em Math. Ann.}, 290(2):209--245, 1991.

\bibitem{BMS05}
A.~I. Bobenko, C.~Mercat, and Y.~B. Suris.
\newblock \href{http://www.arxiv.org/math.DG/0402097}{Linear and nonlinear}
  theories of discrete analytic functions. {I}ntegrable structure and
  isomonodromic {G}reen's function.
\newblock {\em J. Reine Angew. Math.}, 583:117--161, 2005.

\bibitem{BS07}
A.~I. Bobenko and B.~A. Springborn.
\newblock \href{http://arXiv.org/abs/math.DG/0503219}{A discrete
  {L}aplace-{B}eltrami operator} for simplicial surfaces.
\newblock {\em Discrete Comput. Geom.}, 38(4):740--756, 2007.

\bibitem{BS01}
P.~Buser and R.~Silhol.
\newblock Geodesics, periods, and equations of real hyperelliptic curves.
\newblock {\em Duke Math. J.}, 108(2):211--250, 2001.

\bibitem{BS05}
P.~Buser and R.~Silhol.
\newblock Some remarks on the uniformizing function in genus 2.
\newblock {\em Geom. Dedicata}, 115:121--133, 2005.

\bibitem{CSMcC02}
R.~Costa-Santos and B.~M. McCoy.
\newblock \href{http://arXiv.org/abs/hep-th/0109167}{Dimers and the critical
  {I}sing model} on lattices of genus {$>1$}.
\newblock {\em Nuclear Phys. B}, 623(3):439--473, 2002.

\bibitem{CSMcC03}
R.~Costa-Santos and B.~M. McCoy.
\newblock Finite size corrections for the {I}sing model on higher genus
  triangular lattices.
\newblock {\em J. Statist. Phys.}, 112(5-6):889--920, 2003.

\bibitem{DHBvHS04}
B.~Deconinck, M.~Heil, A.~Bobenko, M.~van Hoeij, and M.~Schmies.
\newblock {Computing {R}iemann theta functions}.
\newblock {\em Math. Comp.}, 73(247):1417--1442, 2004.
\newblock
  \href{http://www.ams.org/mcom/2004-73-247/S0025-5718-03-01609-0}{ams}.

\bibitem{DKT08}
M.~Desbrun, E.~Kanso, and Y.~Tong.
\newblock \href{http://www.geometry.caltech.edu/pubs/DKT05.pdf}{Discrete}
  \href{http://www.geometry.caltech.edu/pubs/DKT05.pdf}{differential}
  \href{http://www.geometry.caltech.edu/pubs/DKT05.pdf}{forms}
  \href{http://www.geometry.caltech.edu/pubs/DKT05.pdf}{for computational
  modeling}.
\newblock In A.~I. Bobenko, P.~Schröder, J.~M. Sullivan, and G.~M. Ziegler,
  editors, {\em Discrete Differential Geometry}, volume~38 of {\em Oberwolfach
  Seminars}, pages 287--323. Birkhaüser, 2008.

\bibitem{DMA02}
M.~Desbrun, M.~Meyer, and P.~Alliez.
\newblock \href{http://citeseer.ist.psu.edu/desbrun02intrinsic.html}{Intrinsic
  Parameterizations} of surface meshes.
\newblock In {\em Computer Graphics Forum}, 21, pages 209--218. Proc.
  Eurographics, 2002.

\bibitem{D68}
R.~J. Duffin.
\newblock Potential theory on a rhombic lattice.
\newblock {\em J. Combinatorial Theory}, 5:258--272, 1968.

\bibitem{F44}
J.~Ferrand.
\newblock Fonctions pr\'eharmoniques et fonctions pr\'eholomorphes.
\newblock {\em Bull. Sci. Math. (2)}, 68:152--180, 1944.

\bibitem{FSSB07}
M.~Fisher, B.~Springborn, P.~Schr{\"o}der, and A.~I. Bobenko.
\newblock An algorithm for the construction of intrinsic {D}elaunay
  triangulations with applications to digital geometry processing.
\newblock {\em Computing}, 81(2-3):199--213, 2007.

\bibitem{JWYG04}
M.~Jin, Y.~Wang, S.-T. Yau, and X.~Gu.
\newblock Optimal global conformal surface parameterization.
\newblock In {\em \href{http://www.ic.sunysb.edu/Stu/mijin/vis04.pdf}{VIS '04}:
  Proceedings of the conference on Visualization '04}, 2004.

\bibitem{KSS06}
L.~Kharevych, B.~Springborn, and P.~Schr\"{o}der.
\newblock \href{http://citeseer.ist.psu.edu/kharevych05discrete.html}{Discrete
  conformal} mappings via circle patterns.
\newblock {\em ACM Trans. Graph.}, 25(2):412--438, 2006.

\bibitem{Law70}
H.~B. Lawson, Jr.
\newblock Complete minimal surfaces in {$S\sp{3}$}.
\newblock {\em Ann. of Math. (2)}, 92:335--374, 1970.

\bibitem{M01}
C.~Mercat.
\newblock Discrete {R}iemann surfaces and the {I}sing model.
\newblock {\em Comm. Math. Phys.}, 218(1):177--216, 2001.

\bibitem{M04}
C.~Mercat.
\newblock \href{http://fr.arxiv.org/abs/math-ph/0210016}{Exponentials form a
  basis} of discrete holomorphic functions on a compact.
\newblock {\em Bull. Soc. Math. France}, 132(2):305--326, 2004.

\bibitem{M07}
C.~Mercat.
\newblock Discrete {R}iemann surfaces.
\newblock In A.~Papadopoulos, editor, {\em Handbook of Teichm\"uller Theory,
  vol. I}, volume~11 of {\em IRMA Lect. Math. Theor. Phys.}, pages 541--575.
  Eur. Math. Soc., Z\"urich, 2007.
\newblock \href{http://fr.arxiv.org/abs/0802.1612}{arXiv:0802.1612}.

\bibitem{PP93}
U.~Pinkall and K.~Polthier.
\newblock \href{http://citeseer.ist.psu.edu/518838.html}{Computing discrete
  minimal surfaces} and their conjugates.
\newblock {\em Experiment. Math.}, 2(1):15--36, 1993.

\bibitem{RGA97}
R.~E. Rodr{\'{\i}}guez and V.~Gonz{\'a}lez-Aguilera.
\newblock Fermat's quartic curve, {K}lein's curve and the tetrahedron.
\newblock In {\em Extremal Riemann surfaces (San Francisco, CA, 1995)}, volume
  201 of {\em Contemp. Math.}, pages 43--62. Amer. Math. Soc., Providence, RI,
  1997.

\bibitem{S}
R.~Silhol.
\newblock Period matrices and the {S}chottky problem.
\newblock In {\em Topics on Riemann surfaces and Fuchsian groups (Madrid,
  1998)}, volume 287 of {\em London Math. Soc. Lecture Note Ser.}, pages
  155--163. Cambridge Univ. Press, Cambridge, 2001.

\bibitem{S06}
R.~Silhol.
\newblock Genus 2 translation surfaces with an order 4 automorphism.
\newblock In {\em The geometry of Riemann surfaces and abelian varieties},
  volume 397 of {\em Contemp. Math.}, pages 207--213. Amer. Math. Soc.,
  Providence, RI, 2006.

\bibitem{TACSD06}
Y.~Tong, P.~Alliez, D.~Cohen-Steiner, and M.~Desbrun.
\newblock \href{http://www.eg.org/EG/DL/WS/SGP/SGP06/201-210.pdf}{Designing}
  quadrangulations with discrete harmonic forms.
\newblock In A.~Sheffer and K.~Polthier, editors, {\em Symposium on Geometry
  Processing}, pages 201--210, Cagliari, Sardinia, Italy, 2006. Eurographics
  Association.

\bibitem{Wad06}
M.~Wardetzky.
\newblock {\em Discrete Differential Operators on Polyhedral Surfaces -
  Convergence and Approximation}.
\newblock PhD thesis, \href{http://www.diss.fu-berlin.de/2007/663/}{FU Berlin},
  2006.

\bibitem{W86}
H.~C. Wente.
\newblock Counterexample to a conjecture of {H}. {H}opf.
\newblock {\em Pacific J. Math.}, 121(1):193--243, 1986.

\end{thebibliography}


\end{document}